# Preserving Problem of Local Boundedness of Hamiltonian Dynamical Systems by Symplectic Discretization

Dr. Vice-Prof. *Jong U Hwan, Jo Yon Hui*

Natural Science Centre, **Kim Il Sung** University

**Abstract** We have researched the condition for symplectic discretization to preserve local boundedness
for the space of 2-dimensional Hamiltonian dynamical systems in this paper.

The great leader comrade **Kim Jong Il** said as follows.

"**Scientists and technicians must provide satisfactory solutions to the scientific and technical problems arising in making the national economy Juche-orientated, modern and scientifically-based.**"("**Kim Jong Il** Selected Works" Vol. 10, 194p)

In this paper, we have researched preserving problem of a qualitative structure of the discrete dynamical system which is given by discretization of Hamiltonian dynamical system.

Paper [1] has considered higher order symplectic discretization to preserve energy for Hamiltonian dynamical systems and paper [3] has considered symplectic discretization to preserve energy and momentum for Hamiltonian dynamical systems. And paper [5] has researched explicit discretization to preserve area in phase space and Hamiltonian for 2-dimentional Hamiltonian dynamical system whose potential is cubic or quartic polynomial and paper [6] has derived a explicit symplectic discretization to preserve volume for generalized trigonometric Hamiltonian dynamical system. And paper [4] has considered the condition for various explicit symplectic discretizations to preserve boundedness of orbits about Hamiltonian dynamical system of harmonic oscillator where every orbit is bounded closed curve. The other hand, paper [10] has revealed that first integral of hyper integrable Hamiltonian system is not preserved generally by symplectic discretization.

We have revealed the condition for symplectic discretization to preserve local boundedness for the space of 2-dimensional Hamiltonian dynamical systems in this paper.

We denote the set of non-negative real numbers, positive real numbers and non-positive real numbers by $\mathbf{R}_+$, $\mathbf{R}^+$, $\mathbf{R}_-$ respectively and we denote the set of complex numbers that real part is negative, non-positive by $\mathbf{C}^+$, $\mathbf{C}_-$ respectively. $e$ and $d$ mean natural numbers and $r$ means a natural number or $\infty$ in this paper..

**Definition 1** Suppose $f : \mathbf{R}^e \to \mathbf{R}^e$ is $C^r$ map $(r \geq 1)$. We consider the continuous dynamical system $\varphi : \mathbf{R} \times \mathbf{R}^e \to \mathbf{R}^e$ derived from autonomous system

$$\frac{dy}{dt} = f(y) . \tag{1}$$

Then we call the set $O(x) = \{\varphi(t, x): t \in \mathbf{R}\}$ as orbit of $x$ and we call the set $O_+(x) = \{\varphi(t, x): t \in \mathbf{R}_+\}$, $O_-(x) = \{\varphi(t, x): t \in \mathbf{R}_-\}$ as positive semi-orbit or negative semi-orbit of $x$ respectively. If $f(x) = 0$ then $x \in \mathbf{R}^e$ is called by equilibrium point of (1) and we denote the set of all equilibrium points of (1) by $Eqb$ (1).

**Definition 2** Suppose $F: \mathbf{R}^e \times \mathbf{R}^e \times \mathbf{R}^+ \to \mathbf{R}^e$; $(x, y, \tau) \mapsto F(x, y, \tau)$ is $C^r$ map. We consider one-step discretization for (1) whose order is at least 1

$$y_{n+1} = F(y_n, y_{n+1}, \tau) \tag{2}$$

(we denote (2) by $(2)_\tau$, when $\tau$ is fixed).

If $(2)_\tau$ is expressed like $y_{n+1} = G_\tau(y_n)$ with certain map $G_\tau$ which maps pint of $\mathbf{R}^e$ to point of $\mathbf{R}^e$, then we call set the $O_+(x) = \{G_\tau{}^n(x): n \in \mathbf{Z}_+\}$ as positive semi-orbit of $x$ for discrete dynamical system $(2)_\tau$. We call the point $x \in \mathbf{R}^e$ to satisfy $F(x, x, \tau) = x$ as fixed pint of $(2)_\tau$. We denote the set of all fixed points of $(2)_\tau$ by $Fix(2)_\tau$.

If there exists $\tau_0 > 0$ such that for any $\tau \in (0, \tau_0)$ $Fix(2)_\tau = Eqb(1)$ holds, then we speak like that scheme (2) preserve the equilibrium point set of (1) in the extent $(0, \tau_0)$ of step size.

Now we assume that (2) preserves set of equilibrium points of (1) in the extent $(0, \tau_1)$ of step size. Let's $x_0 \in \mathbf{R}^e$ is equilibrium point of (1). We suppose that linearization of (1) at equilibrium point $x_0$ is given by

$$\frac{dy}{dt} = Ay \tag{3}$$

for $e \times e$ matrix $A$. Thus $A$ is the Jacobian matrix of $f$ at $x_0$.

We assume that the scheme (2) applied to dynamical system (3) for any $\tau \in (0, \tau_1)$ is expressed such as

$$y_{n+1} = S(\tau) y_n \tag{4}$$

for an $e \times e$ matrix $S(\tau)$ (we denote (4) by $(4)_\tau$, when $\tau$ is fixed).

Let's fix $\tau \in (0, \tau_1)$ arbitrarily.

**Definition 3** We denote positive semi-orbit of initial point $x \in \mathbf{R}^e$ of continuous dynamical system (3) by $O_A^+(x)$. And we denote positive semi-orbit of initial point $x$ of discrete dynamical system $(4)_\tau$ by $O_{S(\tau)}^+(x)$ or simply $O_S^+(x)$. We call the set $\mathbf{B}_A = \{x \in \mathbf{R}^e \mid O_A^+(x) \text{ is bounded}\}$ as subspace of dynamical system (3) that positive semi-orbit is bounded. And we call the set $\mathbf{B}_{S(\tau)} = \{x \in \mathbf{R}^e \mid O_{S(\tau)}^+(x) \text{ is bounded}\}$ as subspace of dynamical system $(4)_\tau$ that positive semi-orbit is bounded. We denote $\mathbf{B}_{S(\tau)}$ by simply $\mathbf{B}_S$ occasionally. $\mathbf{B}_A$ and $\mathbf{B}_{S(\tau)}$ is a subspace of $\mathbf{R}^e$ respectively.

**Definition 4** We say that discretization (2) about dynamical system (1) preserves local boundedness at equilibrium point $x_0$ in the extent $(0, \tau_0)$ of step size, if there exists certain $\tau_0 \in (0, \tau_1]$ such that for any $\tau \in (0, \tau_0)$,

1) for any $x \in \mathbf{R}^e$ that $O_A^+(x)$ is bounded, $O_{S(\tau)}^+(x)$ is bounded.

2) if there exists unbounded positive semi-orbit of (3), then $\dim \mathbf{B}_A = \dim \mathbf{B}_{S(\tau)}$.

If discretization (2) for dynamical system (1) preserves local boundedness at equilibrium point $x_0$ in the extent $(0, \tau_0)$ of step size and does not preserve local boundedness at equilibrium point $x_0$ for any $\tau > \tau_0$, then we call $\tau_0$ as preserving limitation of local boundedness at equilibrium point $x_0$.

**Definition 5** If there exists $\tau_0 \in (0, \tau_1]$ such that discretization (2) for dynamical system (1) preserves local boundedness at any $x_0 \in Eqb(1)$ in the extent $(0, \tau_0)$ of step size, then we say that discretization (2) preserve local boundedness for dynamical system (1) at equilibrium points or we simply say that (2) preserve local boundedness.

**Definition 6** Suppose $\mathbf{Q} \subset C^r(\mathbf{R}^e, \mathbf{R}^e)$ is a set of certain dynamical systems. If for any dynamical systems belonging in $\mathbf{Q}$, discretization (2) preserves local boundedness, then we say discretization (2) preserves local boundedness with respect to the set of dynamical systems $\mathbf{Q}$.

**Definition 7** We denote the set of all Hamiltonian dynamical systems

$$\frac{dp}{dt} = -\frac{\partial H}{\partial q}(p, q), \quad \frac{dq}{dt} = \frac{\partial H}{\partial p}(p, q) \tag{5}$$

for $C^2$ map $H : \mathbf{R}^{2d} \to \mathbf{R}$ by $\mathbf{H}(d)$.

And we denote the set of all Hamiltonian dynamical systems

$$\frac{dp}{dt} = -\frac{\partial V}{\partial q}(q), \quad \frac{dq}{dt} = \frac{\partial T}{\partial p}(p) \tag{6}$$

for separable Hamiltonian $H(p, q) = T(p) + V(q)$ (where $T, V : \mathbf{R}^d \to \mathbf{R}$ are $C^2$ map) by $\mathbf{SH}(d)$.

We denote the set of all Hamiltonian dynamical systems (where $g : \mathbf{R}^d \to \mathbf{R}$ is $C^1$ map)

$$\frac{dp}{dt} = g(q), \quad \frac{dq}{dt} = p \tag{7}$$

by $\mathbf{NH}(d)$.

Let's consider symplectic scheme whose order is at least 1 for Hamiltonian dynamical system (5)
$$p_{n+1} = F(p_n, q_n, p_{n+1}, q_{n+1}, \tau), \ q_{n+1} = G(p_n, q_n, p_{n+1}, q_{n+1}, \tau) \tag{8}$$
(We denote (8) by $(8)_\tau$, when $\tau$ is fixed).

Here $F : \mathbf{R}^{4d} \times \mathbf{R}^+ \to \mathbf{R}^d$, $G : \mathbf{R}^{4d} \times \mathbf{R}^+ \to \mathbf{R}^d$ are $C^1$ maps depending on the function $\nabla H(p, q) = \left( \frac{\partial H}{\partial p}(p, q), \ \frac{\partial H}{\partial q}(p, q) \right)$.

Symplecticity of scheme (2) means that $\frac{\partial(p_{n+1}, q_{n+1})}{\partial(p_n, q_n)}^T J \frac{\partial(p_{n+1}, q_{n+1})}{\partial(p_n, q_n)} = J$ holds, where $J$ is the $2d \times 2d$ matrix $J = \begin{pmatrix} 0 & E \\ -E & 0 \end{pmatrix}$ and $E$ is $d \times d$ unit matrix([2]).

Let's consider the case $d = 1$. Let linearization of Hamiltonian dynamical system (5) at equilibrium point $(p_0, q_0) \in \mathbf{R}^2$ is given by

$$\begin{pmatrix} \dfrac{dp}{dt} \\ \dfrac{dq}{dt} \end{pmatrix} = A \begin{pmatrix} p \\ q \end{pmatrix} \tag{9}$$

where $A$ is a $2 \times 2$ matrix. And we assume that symplectic scheme (2) applied to (9) is expressed by

$$\begin{pmatrix} p_{n+1} \\ q_{n+1} \end{pmatrix} = S(\tau) \begin{pmatrix} p_n \\ q_n \end{pmatrix} \tag{10}$$

where $S(\tau)$ is a $2 \times 2$ matrix.

**Theorem 1** Let's $d = 1$. If there exists $\tau_0 \in (0, \tau_1]$ such that for any $\tau \in (0, \tau_0)$ following conditions 1)-4) hold, then scheme (2) preserves local boundedness about Hamiltonian dynamical system (1) at equilibrium point $(p_0, q_0)$ in the extent $(0, \tau_0)$ of step size.

1) if $\det A > 0$, then $|Tr\, S(\tau)| < 2$.

2) if $\det A < 0$, then $|Tr\, S(\tau)| > 2$.

3) if $rank\, A = 1$, then $Tr\, S(\tau) = 2$, $rank\,(S(\tau) - E) = 1$.

4) if $rank\, A = 0$, then $Tr\, S(\tau) = 2$, $rank\,(S(\tau) - E) = 0$.

**Proof** $Tr A = a_{11} + a_{22} = 0$ holds because (9) is Hamiltonian equation and $\det S(\tau) = 1$ holds because (8) is symplectic scheme. If we compare the classification of orbits of (9) based on $Tr A = a_{11} + a_{22} = 0$ with the classification of orbits of (10) based on $\det S(\tau) = 1$, then we get the result of theorem. Q.E.D.

**Corollary 1** We assume that $\det A < 0$ holds in (9). We suppose that there exists $0 < \tau_0 < \tau_1$ such that for any $\tau \in (0, \tau_0)$, $|Tr\, S(\tau)| > 2$ holds. Then the numerical solution by symplectic scheme $(2)_\tau$ for $\tau \in (0, \tau_0)$ applied to (9) is stable([7, 9]). Thus the error between real solution and approximate solution $Y_n = y_n - y(n\tau)$, $n \geq 0$ is bounded.

**Proof** We denote the solution of (9) satisfying initial condition $y(0) = y^0$ by $y(t)$. We assume rounding error $\phi$ emerges at every step which is independent with step number $n$ when we calculate numerical solution. And we assume that local discretization error $\varepsilon$ is independent with step number $n$. Then we have

$$y_{n+1} = S(\tau) y_n + \phi, \quad y((n+1)\tau) = S(\tau) y(n\tau) + \tau \varepsilon(\tau).$$

If we subtract upper equation from lower equation, then we get equation for error $Y_n$ such as

$$Y_{n+1} = S(\tau) Y_n + \eta(\tau). \tag{11}$$

Here $\eta(\tau) = \phi - \tau \varepsilon(\tau)$. The solution of difference equation (11) satisfying initial condition $Y_0 = y_0 - y^0$ is given by

$$Y_n = S(\tau)^n [y_0 - y^0 - (I - S(\tau))^{-1} \eta(\tau)] + (I - S(\tau))^{-1} \eta(\tau). \tag{12}$$

Any solution of (9) is bounded from the assumption $\det A < 0$. The solution $Y_n = S(\tau)^n \xi$, $n \geq 0$ of difference equation $Y_{n+1} = S(\tau) Y_n$, $Y_0 = \xi$ is bounded for any $\xi \in \mathbf{R}^2$ under the condition $|Tr\, S(\tau)| > 2$ ($\tau \in (0, \tau_0)$) from the theorem. Therefore $S(\tau)^n [y_0 - y^0 - (I - S(\tau))^{-1} \eta(\tau)]$, $n \geq 0$ is bounded. Thus $\{Y_n, n \geq 0\}$ is bounded. Q.E.D.

**Corollary 2** Symplectic Euler B method ([2, 7])

$$P = p - \frac{\partial V}{\partial q}(Q), \ Q = q + \frac{\partial T}{\partial p}(p)$$

for separable Hamiltonian dynamical system preserves local boundedness with respect to the space **SH**(1) of Hamiltonian dynamical systems.

Moreover if $(p_0, \ q_0)$ is equilibrium point of separable.

Hamiltonian dynamical system (6), then the limitation of preserving local boundedness at equilibrium point $(p_0, \ q_0)$ of (6) is given by

$$\tau_{\max} = \frac{2}{\sqrt{\dfrac{\partial^2 T}{\partial p^2}(p_0)\dfrac{\partial^2 V}{\partial q^2}(q_0)}}, \ (\text{if } \ T''(p_0)V''(q_0) > 0 )$$

$$\tau_{\max} = +\infty , \qquad\qquad (\text{if } \ T''(p_0)V''(q_0) \le 0).$$

**Corollary 3** Explicit Symplectic scheme of order 2 ([8])

$$\tilde{q} = q + \frac{\tau}{2}\nabla T(p), \ \ \tilde{p} = p - \frac{\tau}{2}\nabla V(\tilde{q}), \ \ Q = \tilde{q} + \frac{\tau}{2}\nabla T(\tilde{p}), \ \ P = \tilde{p}$$

for separable Hamiltonian dynamical system preserves local boundedness with respect to the space **SH**(1) of Hamiltonian dynamical systems. Moreover if $(p_0, \ q_0)$ is equilibrium point of separable

Hamiltonian dynamical system (6), then the limitation of preserving local boundedness at equilibrium point $(p_0, q_0)$ of (6) is given by

$$\tau_{\max} = \frac{2}{\sqrt{T''(p_0)V''(q_0)}}, \ \ (\text{if } \ T''(p_0)V''(q_0) > 0 )$$

$$\tau_{\max} = +\infty , \qquad\qquad (\text{if } \ T''(p_0)V''(q_0) \le 0).$$

**Corollary 4** Stömer-Verlet scheme([2])

$$p_{k+1/2} = p_k + \tau g(q_k)/2, \ \ q_{k+1} = q_k + \tau \, p_{k+1/2}, \ \ p_{k+1} = p_{k+1/2} + \tau g(q_{k+1})/2$$

for Hamiltonian dynamical system (7) preserves local boundedness with respect to the space **NH**(1) of Hamiltonian dynamical systems. Moreover if $(p_0, \ q_0)$ is equilibrium point of separable Hamiltonian dynamical system (7), then the limitation of preserving local boundedness at equilibrium point $(p_0, q_0)$ of (7) is given by

$$\tau_{\max} = \frac{2}{\sqrt{-g'(q_0)}}, \ \ (\text{if } \ g'(q_0) < 0 )$$

$$\tau_{\max} = +\infty , \qquad\qquad (\text{if } \ g'(q_0) \ge 0).$$

**Corollary 5** Implicit Euler scheme ([2])

$$P = p - \tau\frac{\partial H}{\partial q}\left(\frac{p+P}{2}, \ \frac{q+Q}{2}\right), \ \ Q = q + \tau\frac{\partial H}{\partial p}\left(\frac{p+P}{2}, \ \frac{q+Q}{2}\right)$$

for Hamiltonian dynamical system (5) preserves local boundedness about the space **H**(1) of Hamiltonian dynamical systems. Moreover if $(p_0, \ q_0)$ is equilibrium point of separable Hamiltonian dynamical system (5), then the limitation of preserving local boundedness at equilibrium point $(p_0, q_0)$ of (5) is given by

$$\tau_{\max} = \frac{2}{\sqrt{-H_0}}, \quad (\text{if } H_0 < 0)$$

$$\tau_{\max} = +\infty, \qquad (\text{if } H_0 \geq 0)$$

where $H_0 = H_{pp}(p_0, q_0)H_{qq}(p_0, q_0) - H_{pq}(p_0, q_0)^2$.

**Refference**